\documentclass[12pt, a4paper]{article}
\usepackage{natbib}
\usepackage{float}
\usepackage{multirow}
\usepackage{amssymb,amsmath, amsfonts}
\usepackage{lscape}
\usepackage{enumerate}
\usepackage{amsthm}
\usepackage[T1]{fontenc}
\usepackage[utf8]{inputenc}
\usepackage{lmodern}
\usepackage[english]{babel}
\usepackage{csquotes}
\usepackage{mathtools}

\usepackage{xcolor}
\numberwithin{equation}{section}
\usepackage{graphicx}
\usepackage{caption}
\usepackage{subcaption}
\graphicspath{{density.jpeg},{distribution.jpeg},{3mp1mfig,jpeg},{38fig.jpeg},{225fig.jpeg}}

\newtheorem{proposition}{Proposition}[section]

\newtheorem{remark}{Remark}[section]

\usepackage[figuresright]{rotating}
\begin{document}

\begin{center}
{\Large\bf A new class of copulas having dependence range larger than FGM-type copulas} \\
\vspace{0.3in}

{ \bf Swaroop Georgy Zachariah$^{a}$,  \bf Mohd. Arshad$^{a}$\footnote{Corresponding author. E-mail
addresses: ~arshad@iiti.ac.in(M. Arshad), ~ashokiitb09@gmail.com (Ashok Kumar Pathak), ~swaroopgeorgy@gmail.com
(Swaroop Georgy Zachariah).} and \bf Ashok Kumar Pathak$^{b}$ }
 \\
$^a$ Department of Mathematics, Indian Institute of
Technology Indore, India \\
 $^b$ Department of Mathematics and Statistics,
Central University of Punjab, Bathinda,
India
\end{center}

\noindent------------------------------------------------------------------------------------------------------------------------------
\begin{abstract}
We propose a new bivariate symmetric copula with positive and negative dependence properties. The main features of the proposed copula are its simple mathematical structure, wider dependence range compared to FGM copula and its generalizations, and no lower and upper tail dependence. The maximum range of Spearman's Rho of the proposed copula is $[-0.5866,0.5866]$, which improves the dependence range of the FGM copula and its various generalizations. A new bivariate Rayleigh distribution is developed using the proposed copula, and some statistical properties have been studied. A real data set is analyzed to illustrate the proposed bivariate distribution's relevance in practical contexts.
\end{abstract}
\vskip 2mm

 \noindent \emph{Keywords}: Copula, Kendall's Tau, Spearman's Rho, Rayleigh Distribution, Bivariate Distribution.\\

\noindent \emph{Mathematics Subject Classification (2020)}: 62H05, 62H10. \\
\noindent------------------------------------------------------------------------------------------------------------------------------


\section{Introduction}
Copula plays a significant role in the field of statistics, finance, engineering and medical sciences for modelling dependent data sets. If we have a family of copulas, we naturally have a collection of multivariate distributions with whatever marginal distributions we desire. This feature of the copula is useful in every branch of study where dependence modelling and simulation are an integral part. \cite{sklar1959fonctions} proved that if $X$ and $Y$ are two random variables having marginal cumulative distributions $G(x)$ and $H(y)$ respectively, then there exists a function $C$ such that 
$F(x,y)=P(X\leq x, Y\leq y)=C(G(x),H(y)), \ \forall(x,y)\in \mathbb{R}^2$. If $G(\cdot)$ and $H(\cdot)$ are continuous, then $C$ is unique: otherwise $C$ is uniquely determined on Range$(F)\times$Range$(G)$. Thus, every joint cumulative distribution function can be expressed as a function of the cumulative distributions of corresponding marginals via copula. Moreover, a bivariate function $C:\mathcal{I}^2\rightarrow\mathcal{I}$ is said to be a bivariate copula if it satisfies the following conditions:
\begin{align}
	&C(u,0)=C(0,v)=0 \ \ \text{and} \ \ C(u,1)=u,C(1,v)=v, \ \forall u,v \in \mathcal{I}\label{boundarycond}\\
	&C(u_2,v_2)-C(u_1,v_2)-C(u_2,v_1)+C(u_1,v_1)\geq 0, \ \forall u_1<u_2,v_1<v_2 \in \mathcal{I} \label{2inc}
\end{align}
where $\mathcal{I}=[0,1]$ and $\mathcal{I}^2=[0,1] \times [0,1]$. Note that Eq. (\ref{2inc}) is called the 2-increasing property. For more details about copulas, see \cite{nelsen2007introduction}, \cite{durante2016principles}, \cite{mai2017simulating}, \citet{hofert2018elements}, and references therein. 

\par In recent years, many researchers have used copula functions to construct multivariate distributions. For instances, \cite{achcar2015bivariate} studied bivariate generalized exponential distribution via Farlie-Gumbel-Morgenstern (FGM) copula. \cite{kundu2017bivariate} discussed bivariate Birnbaum-Saunders using Gaussian copula. \cite{abd2017inference} studied bivariate generalized exponential distribution from FGM and Plackett copula functions and discussed various estimation procedures. \cite{el2018fgm} used FGM copula for constructing a bivariate Weibull distribution. Recently, \cite{almetwally2020bivariate} developed bivariate Fr{\'e}chet distribution using the FGM copula.
\par In literature, a wide variety of copulas are available; of them, FGM copula received much attention due to its simple mathematical structure and exhibited positive and negative dependence (see \cite{farlie1960performance}). Spearman's Rho, denoted by $\rho_c$, is one of the dependence measures used for measuring the dependence structure captured by a copula. The range of values of $\rho_c$ lies between -1 and 1. Positive values of $\rho_c$ indicate positive dependence and negative values for negative dependence. For the FGM copula, the range of values of Spearman's Rho is very low, i.e., $\rho_c\in[-0.33,0.33]$ (see \cite{farlie1960performance}). So, the FGM copula is unsuitable for modelling data with a high dependence structure. Many researchers attempted to propose FGM-type copula for improving the correlation coefficient. \cite{huang1999modifications} proposed two extended FGM copulas, having $\rho_c\in[-0.33,0.375]$ and $\rho_c\in[-0.33,0.391]$ respectively. \cite{bairamov2002dependence} and \cite{bekrizadeh2015pak} also extended FGM copula with $\rho_c\in[-0.48,0.502]$ and $\rho_c\in[-0.50,0.43]$ respectively. Recently, \cite{chesneau2022note} proposed a polynomial-sine copula exhibiting positive as well as negative dependence with $\rho_c\in[-0.4927,0.4927]$. 
In most of these works, we can see that parameters are added to improve Spearman's correlation coefficient range, resulting in a mathematically complex structure and computationally more expensive for estimating unknown parameters. To overcome these drawbacks, we propose a simple bivariate copula without adding any parameter to existing copulas. The proposed copula  improves the dependence range of Spearman's correlation of various FGM-type copulas reported in the literature.
\par The present paper is organized as follows: In Sections \ref{2} and \ref{3}, we will introduce a new bivariate copula and study some dependency measures of the proposed copula. In Section \ref{4} a new bivariate Rayleigh distribution is derived from the proposed copula. We also derive expressions for conditional distribution and product moments. A real data is also analyzed using the proposed bivariate Rayleigh distribution as an application.

\section{New Bivariate Copula}\label{2}
Consider the following bivariate function 
\begin{equation}\label{copula}
	C(u,v;\delta,\alpha)=uv+\delta\left(1-e^{\alpha(u-u^2)}\right)\left(1-e^{\alpha(v-v^2)}\right), \quad (u,v) \in \mathcal{I}^2,
\end{equation}
where $\alpha$ and $\delta$ are real valued parameters. Assume that the parameter $\delta$ depends on $\alpha$, and $\alpha$ takes arbitrary value in $\mathbb{R}$. The bivariate function, defined in Eq. (\ref{copula}), satisfied the boundary conditions of a bivariate copula given in Eq. (\ref{boundarycond}). But we need to find the range of parameter $\delta$ for which this function is a valid bivariate copula, i.e., the function in Eq. (\ref{copula}) satisfy the $2$-increasing property given in Eq. (\ref{2inc}). \cite{kim2011generalized} proved that $2$-increasing property in an absolutely continuous copula is equivalent to the condition that copula density $c(u,v;\delta,\alpha)$ is non-negative, i.e., 
\begin{equation}\label{cd}
	c(u,v;\delta,\alpha)=\dfrac{\partial^2C(u,v)}{\partial u\partial v}=1+\delta \alpha^2g(u)g(v)\geq 0,
\end{equation}
where $g(t)=\left(1-2t\right)e^{\alpha(t-t^2)}, t\in\mathcal{I}$. Clearly, non-negativity of copula density in Eq. (\ref{cd}) depends on the behaviour of function $g(\cdot)$ and the value of the parameter $\delta$.   For finding the feasible range of $\delta$, we will divide the domain $\mathcal{I}^2$ of $(u,v)$ into four quadrants as:
\begin{equation*}
	\begin{aligned}
	R_1 &=\left\{(u,v)\in\mathcal{I}^2:0\leq u\leq \frac{1}{2},0\leq v\leq \frac{1}{2} \right\},
		\\
	R_2 &=\left\{(u,v)\in\mathcal{I}^2:\frac{1}{2}< u\leq 1,0\leq v\leq \frac{1}{2} \right\},
		\\
	R_3 &=\left\{(u,v)\in\mathcal{I}^2:\frac{1}{2}< u\leq 1, \frac{1}{2}< v\leq 1 \right\},
	\\
	R_4&=\left\{(u,v)\in\mathcal{I}^2:0\leq u\leq \frac{1}{2},\frac{1}{2}< v\leq 1 \right\}.
	\end{aligned}
\qquad\qquad \vcenter{\hbox{\includegraphics[width=4.5cm,height=4.5cm]{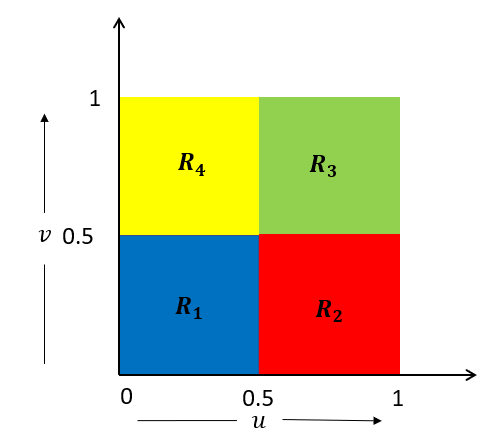}}}
\end{equation*}
Since the sign of $g(t)$ is positive if $t \in \left(0, \frac{1}{2}\right)$, and negative if $t \in \left(\frac{1}{2}, 1\right)$, it follows that the product $g(u)g(v)$ is positive on $R_1 \cup R_3$ and negative on $R_2 \cup R_4$. Thus, the copula density $c(u,v;\delta,\alpha)$ is non-negative if 
\begin{equation}\label{upper}
	\delta\geq \frac{-1}{\alpha^2 g(u)g(v)},  \quad (u,v) \in R_1 \cup R_3, \nonumber
\end{equation} 
and 
\begin{equation}\label{lower}
\delta\leq \frac{-1}{\alpha^2 g(u)g(v)},  \quad (u,v) \in R_2 \cup R_4. \nonumber
\end{equation}
Therefore, the copula density $c(u,v;\delta,\alpha)$ is non-negative if 
\begin{equation}\label{rangeD}
\frac{-1}{\alpha^2 \underset{(u,v) \in R_1 \cup R_3}{\sup} \{g(u)g(v)\}} \leq 	\delta \leq \frac{-1}{\alpha^2 \underset{(u,v) \in R_2 \cup R_4}{\inf} \{g(u)g(v)\}}. 
\end{equation}
Since the behaviour of function $g(\cdot)$ depends on $\alpha \in \mathbb{R}$, we will consider the following three cases: \\
{\bf{Case I:}} When $\alpha=0$. \\
In this case, the copula density $c(u,v;\delta,0)=1$, which is non-negative, and hence the $2$-increasing property holds for arbitrary value of $\delta$. Moreover, the proposed copula (\ref{copula}) reduces to the product copula. Note that the product copula is a well known copula, which corresponds to the independence of two random variables. \\

\begin{figure}[ht]
	\centering
	\subfloat[$\alpha=-3,\delta=-0.1$]{\includegraphics[height=5cm, width=5cm]{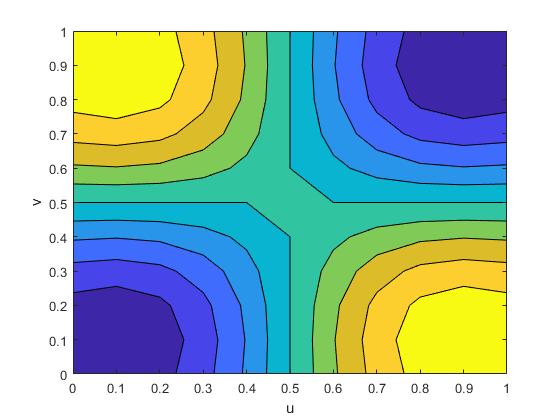}}  	 
	\subfloat[$\alpha=2,\delta=0.2$]{\includegraphics[height=5cm, width=5cm]{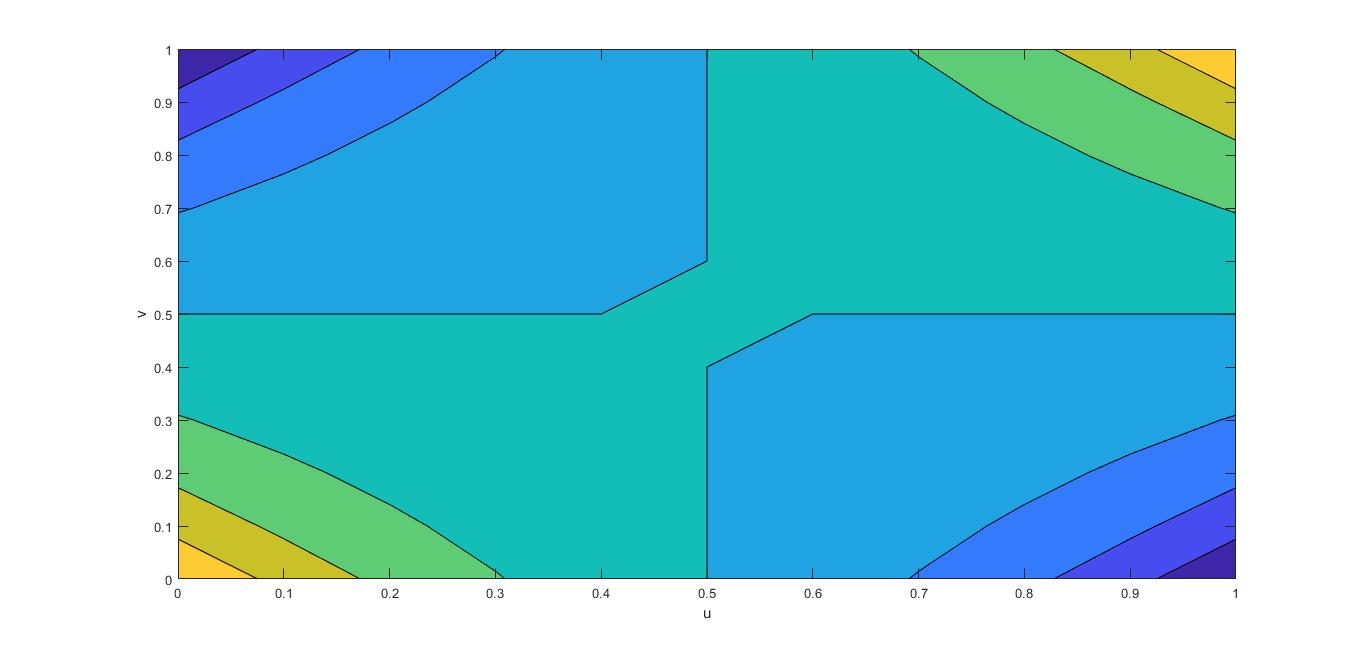}}            
	\subfloat[$\alpha=3.8,\delta=-0.3$]{\includegraphics[height=5cm, width=5cm]{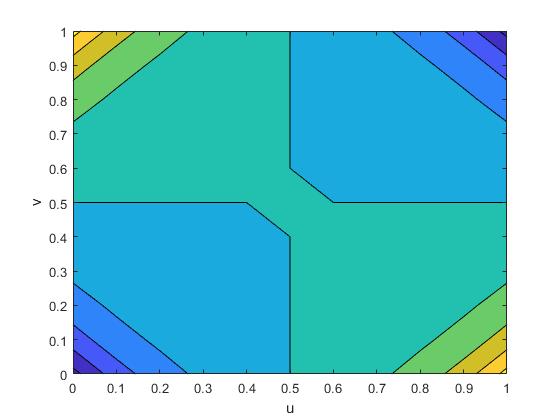}} 
	\caption{Contour plots of copula density $c(u,v)$ for various values of $\alpha$ and $\delta$. 
	}\label{cont}
\end{figure}

\noindent{\bf{Case II:}} When $\alpha\leq2, \alpha\neq0$.\\
In this case, $g(t)$ is a decreasing function on $\mathcal{I}$ with $g(0)=1, g\left(\frac{1}{2}\right)=0$ and $g(1)=-1$. It follows that $g(t)$ takes values in $[0,1]$ for $t \in \left[0, \frac{1}{2}\right]$, and takes values in $\left[-1,0\right)$ for $t \in \left(\frac{1}{2}, 1\right]$. Thus, the product function $g(u)g(v)$ is bounded above by $1$ on  $R_1 \cup R_3$, and the upper bound $1$ is attended at $(u,v) \in \{(0,0), (1,1)\}$. Therefore, $\underset{(u,v) \in R_1 \cup R_3}{\sup} \{g(u)g(v)\}=1$. Further, for $(u,v) \in R_2$, $g(u)$ takes values in $[-1,0)$ and $g(v)$ takes values in $[0,1]$. This implies that the product function $g(u)g(v)$ is bounded below by $-1$ on $R_2$, and the lower bound $-1$ is attended at $u=1,v=0$. Similarly, for $(u,v) \in R_4$, the product function $g(u)g(v)$ is bounded below by $-1$, which is attended at $u=0,v=1$. Therefore, $\underset{(u,v) \in R_2 \cup R_4}{\inf} \{g(u)g(v)\}=-1$. Now, using these values in inequality (\ref{rangeD}), we get the feasible range of the parameter $\delta$ as 
\begin{equation}\label{rangeD1}
	\frac{-1}{\alpha^2} \leq 	\delta \leq \frac{1}{\alpha^2}. 
\end{equation}  
\noindent{\bf{Case III:}} When $\alpha>2$.\\
Let $r_1=\frac{1}{2}-\frac{1}{\sqrt{2 \alpha}}$ and let $r_2=\frac{1}{2}+\frac{1}{\sqrt{2 \alpha}}$. Clearly, $0<r_1<\frac{1}{2}<r_2<1$. It can be observed that $g(t)$ increases on $t \in [0,r_1]$, decreases on $t \in (r_1,r_2)$, and increases on $t \in [r_2,1]$. Also, $g(t)$ takes positive values on $t \in \left[0, \frac{1}{2}\right)$ and negative values on $t \in \left(\frac{1}{2}, 1\right]$, with $g\left(\frac{1}{2}\right)=0$. Moreover, $g(t)$ is maximum at $t=r_1$ with maximum value $g(r_1)=\sqrt{\frac{2}{\alpha}} \exp\left\{\frac{\alpha}{4}-\frac{1}{2}\right\}$, and $g(t)$ is minimum at $t=r_2$ with minimum value $g(r_2)=-\sqrt{\frac{2}{\alpha}} \exp\left\{\frac{\alpha}{4}-\frac{1}{2}\right\}$. Since the functions $g(u)$ and $g(v)$ are positive on the quadrant $R_1$ and takes maximum at $u=r_1, v=r_1$, it follows that the product function $g(u)g(v)$ has maximum value $\left[g(r_1)\right]^2=\frac{2}{\alpha} \exp\left\{\frac{\alpha}{2}-1\right\}$, on $(u,v) \in R_1$. Since the functions $g(u)$ and $g(v)$ are negative on $R_3$ and takes minimum at $u=r_2, v=r_2$, it follows that the product function $g(u)g(v)$ is positive and has maximum value $\left[g(r_2)\right]^2=\frac{2}{\alpha} \exp\left\{\frac{\alpha}{2}-1\right\}$, on $(u,v) \in R_3$. Therefore, $\underset{(u,v) \in R_1 \cup R_3}{\sup} \{g(u)g(v)\}=\frac{2}{\alpha} \exp\left\{\frac{\alpha}{2}-1\right\}$. Similarly, we have found that the infimum of the product function $g(u)g(v)$ on the quadrant $R_2 \cup R_4$ is equal to $g(r_1)g(r_2)=-\frac{2}{\alpha} \exp\left\{\frac{\alpha}{2}-1\right\}$. Now, using these values in Eq. (\ref{rangeD}), we get the feasible range of $\delta$ as
\begin{equation}\label{rangeD2}
-\frac{1}{2\alpha}\exp\left\{1-\frac{\alpha}{2}\right\}  \leq \delta  \leq\frac{1}{2\alpha}\exp\left\{1-\frac{\alpha}{2}\right\}.
\end{equation}  
Thus, using the feasible range of the parameter $\delta$ given in Eq. (\ref{rangeD1}) and Eq. (\ref{rangeD2}), we propose the following bivariate copula  
\begin{equation}\label{Copula_prop}
	C(u,v;\delta,\alpha)=uv+\delta\left(1-e^{\alpha(u-u^2)}\right)\left(1-e^{\alpha(v-v^2)}\right), \quad (u,v) \in \mathcal{I}^2, 
\end{equation}
where $\alpha$ is a real valued parameter, and $\lvert\delta\rvert\leq \delta ^{\star}(\alpha)$. Here, 
\begin{align*}
	\delta ^{\star}(\alpha)=&\left \{ \begin{array}{cl} \   \frac{1}{\alpha^2},& \mbox{if} \ \alpha\in\left(-\infty,2\right]\setminus\left\{0\right\}\\
		\\
		\frac{1}{2\alpha}\exp\left\{1-\frac{\alpha}{2}\right\}, & \mbox{if} \ \alpha>2.
	\end{array} \right.
\end{align*}
Recall that our proposed copula reduced to the product copula when $\alpha=0$. Also, the copula density $c(u,v)=\dfrac{\partial^2 C(u,v)}{\partial u\partial v}$ of the proposed copula is given by 
\begin{equation}\label{copulad}
	c(u,v)=1+\alpha^2\delta\left(1-2u\right)\left(1-2v\right)\exp\left\{\alpha\left(u-u^2+v-v^2\right)\right\}, \quad (u,v) \in\mathcal{I}^2, 
\end{equation}
where $\alpha\in\mathbb{R}$ and $|\delta|\leq\delta^{\star}(\alpha)$. The contour plots of the copula density are shown in Figure \ref{cont} for different choices of the parameters $\alpha$ and $\delta$.


%


\section{Measures of Dependence}\label{3}
\par Copula functions have been widely used for modelling dependence between random variables since they allow the separation of the dependence effect from the effects of the marginal distributions. The dependence between random variables can be measured by several well-known measures of dependence such as Spearman's Rho, Gini's Gamma Coefficient, Kendall's tau, Spearman's footrule coefficient (see \cite{nelsen2007introduction}) and Blest's measure of rank correlation (see \citealp{genest2003blest}). These measures are usually based on the ranks of observations rather than their actual values. The range of these dependence measures lies between $-1$ and $1$. A negative value of dependence measure corresponds to negative dependence, zero for independence and a positive value for positive dependence. Most of these measures remain unchanged under strictly increasing transformations of random variables. Since the copula functions of a pair of random variables $X$ and $Y$ are invariant under strictly increasing transformations of $X$ and $Y$, these dependence measures are expressible in terms of the copulas. Let $C(u,v)$ be the copula function associated with the random variables $X$ and $Y$. The following are important measures of dependence in terms of copula function: 
\begin{itemize}
\item {\bf  Spearman's Rho}  
\begin{equation}
	\rho_c=12\int_{0}^{1}\int_{0}^{1}C(u,v)dudv-3. \nonumber \label{spearman}
\end{equation}

\item {\bf Gini's Gamma Coefficient}  
\begin{equation}
	\gamma_c=4\left\{\int_{0}^{1}C(u,1-u)du-\int_{0}^{1}\left(u-C(u,u)\right)du\right\}.
	\nonumber 
\end{equation} 
	\item {\bf Kendall's Tau}  
\begin{equation}
	\tau_c=4\int_{0}^{1}\int_{0}^{1}C(u,v)dC(u,v)-1. \nonumber 
\end{equation}
\item {\bf Blest's Measure of Rank Correlation}

\begin{equation}
	\eta_c=24\int_{0}^{1}\int_{0}^{1}\left(1-u\right)C(u,v)dudv-2. \nonumber 
\end{equation}

\item {\bf Spearman's footrule}
\begin{equation}
	\phi_c=6\int_{0}^{1}C(u,u)du-2. \nonumber
\end{equation}
\end{itemize}
Now, we will provide the expressions of the various measures of dependence for the proposed copula function (\ref{Copula_prop}).  
\begin{proposition}
For the copula defined in Eq. (\ref{Copula_prop}), the Spearman's Rho $\rho_C$, and Gini's Gamma Coefficient $\gamma_C$ are given by
	\begin{align*}
	 \rho_C=&\left \{ \begin{array}{cl} \  \displaystyle 12\delta\left[1-\sqrt{\frac{\pi}{|\alpha|}}e^{\alpha/4}erfi\left\{\frac{\sqrt{|\alpha|}}{2}\right\}\right]^2,& \mbox{if} \ \alpha<0\\
			\\
			12\delta\left[1-\sqrt{\frac{\pi}{\alpha}}e^{\alpha/4}erf\left\{\frac{\sqrt{\alpha}}{2}\right\}\right]^2, & \mbox{if} \ \alpha \geq0,\end{array} \right.\\
		\gamma_C=&\left \{ \begin{array}{cl} \  \displaystyle 8\delta\left[1-2\sqrt{\frac{\pi}{|\alpha|}}e^{\alpha/4}erfi\left\{\frac{\sqrt{|\alpha|}}{2}\right\}+\sqrt{\frac{\pi}{2|\alpha|}}e^{\alpha/2}erfi\left\{\sqrt{\frac{|\alpha|}{2}}\right\}\right],& \mbox{if} \ \alpha<0\\
			\\
		8\delta\left[1-2\sqrt{\frac{\pi}{\alpha}}e^{\alpha/4}erf\left\{\frac{\sqrt{\alpha}}{2}\right\}+\sqrt{\frac{\pi}{2\alpha}}e^{\alpha/2}erf\left\{\sqrt{\frac{\alpha}{2}}\right\}\right], & \mbox{if} \ \alpha \geq0,\end{array} \right.\\
	\end{align*}
	where $erf(t)=\frac{2}{\sqrt{\pi}}\int_0^{t}e^{-z^2}dz$ denotes the error function (see \cite{abramowitz1972handbook}) and $erfi(t)=\frac{2}{\sqrt{\pi}}\int_0^{t}e^{z^2}dz$ denotes the imaginary error function (see \cite{marcinowski2020using}).
\end{proposition}
\begin{remark}
	It can be verified that the expressions of other measures of dependence satisfied the following relation under the copula given in Eq. (\ref{Copula_prop}). 
	\[\eta_c=\rho_c=\frac{3}{2}\tau_c, \quad {\text{and}} \quad \phi_c=\frac{3}{4}\gamma_c. \]
\end{remark}
Table \ref{depende} presents the numerical values of Spearman's rho and Gini's gamma coefficient of the new copula for different values of the copula parameter $\alpha$. Since the new copula is symmetric, we have shown only the upper boundary values of the dependence measures in Table \ref{depende}. The lower boundary value is the negative of upper boundary value. It is observed from the Table \ref{depende} that Spearman's rho $\rho_c\in[-0.5866,0.5866]$ when $\alpha=3.8$, thereby extending the range of the dependence measure Spearman's rho over the popular FGM copula and its various generalizations. 
\begin{table}[h]
	\caption{Sperman's Rho and Gini's Gamma Coefficient for various values of $\alpha$} 
	\centering 
		\begin{tabular}{cccccccc}
			\hline
			$\alpha$ & $\delta_{upper}$ & $\rho_{upper}$ & $\gamma_{upper}$ & $\alpha$ & $\delta_{upper}$ & $\rho_{upper}$ & $\gamma_{upper}$ \\ \hline
			-3   & 0.1111  & 0.1899 & 0.1463 & 1.2 & 0.6944 & 0.4264 & 0.3473 \\
			-2.7 & 0.1372  & 0.2003 & 0.1547 & 1.5 & 0.4444 & 0.4544 & 0.3718 \\
			-2.4 & 0.1736  & 0.2113 & 0.1638 & 1.8 & 0.3086 & 0.4845 & 0.3984 \\
			-2.1 & 0.2268  & 0.2231 & 0.1736 & 2   & 0.25   & 0.506  & 0.4174 \\
			-1.8 & 0.3086  & 0.2358 & 0.1841 & 2.3 & 0.1871 & 0.5348 & 0.4434 \\
			-1.5 & 0.4444  & 0.2493 & 0.1954 & 2.6 & 0.1425 & 0.5561 & 0.4634 \\
			-1.2 & 0.6944  & 0.2638 & 0.2076 & 2.9 & 0.1099 & 0.571  & 0.4783 \\
			-0.9 & 1.2346  & 0.2794 & 0.2207 & 3.2 & 0.0858 & 0.5805 & 0.4888 \\
			-0.6 & 2.7778  & 0.2961 & 0.2349 & 3.5 & 0.0675 & 0.5855 & 0.4956 \\
			-0.3 & 11.1111 & 0.314  & 0.2502 & 3.8 & 0.0535 & {\bf 0.5866} & 0.4992 \\
			0    & 0       & 0      & 0      & 4.1 & 0.0427 & 0.5845 & {\bf 0.5002} \\
			0.3  & 11.1111 & 0.3541 & 0.2845 & 4.4 & 0.0342 & 0.5798 & 0.499  \\
			0.6  & 2.7778  & 0.3764 & 0.3038 & 4.7 & 0.0276 & 0.5729 & 0.4959 \\
			0.9  & 1.2346  & 0.4005 & 0.3247 & 5   & 0.0223 & 0.5643 & 0.4912 \\ \hline
	\end{tabular}\label{depende}
\end{table}
\par In order to continue our discussion on dependence between random variables, there are some more dependence properties available in the literature. For example, quadrant dependence, totally positive of order 2 ($TP_2$), and tail dependence coefficient. For more detailed discussion on these properties, one can see \cite{lehmann1966some}, \cite{barlow1975statistical}, \cite{drouet2001correlation}, \cite{nelsen2007introduction}, \cite{lai2009continuous} and \cite{bhuyan2020bivariate}. Now, we will discuss these dependence properties under the copula given in Eq. (\ref{Copula_prop}).     
\subsection{Quadrant Dependence} 
\par The random vector $(X, Y)$ is said to be positively (negatively) quadrant dependence if 
$$P(X\leq x,Y\leq y)\geq (\leq) \ P(X\leq x)P(Y \leq y), \quad \forall \;  (x,y)\in\mathbb{R}^2.$$
Similarly, we call copula $C(u,v)$ is positively (negatively) quadrant dependent if 
$$C(u,v)\geq (\leq) \ uv, \quad \forall (u,v)\in\mathcal{I}^2.$$
It can be verified that the quadrant dependence of the copula (\ref{Copula_prop}) depends only on the copula parameter $\delta$. The proof is straightforward, so omitted. Thus, we have the following result. 
\begin{proposition}
	The copula (\ref{Copula_prop}) has positive (negative) quadrant dependence if $\delta\geq0$ \ ($\delta\leq0$).
\end{proposition}
\subsection{Totally Positive of Order 2 ($TP_2$)} 
A bivariate function $g(x,y)$ is said to be totally positive of order 2 ($TP_2$) if 
$$g(x_1,y_1)g(x_2,y_2)-
g(x_2,y_1) g(x_2,y_1)\geq0, \quad \forall \ x_1<x_2, y_1<y_2.$$
The $TP_2$ property is a stronger concept of dependence. If a copula density $c(u,v)$ possesses $TP_2$ property, then the associated copula $C(u,v)$ has stochastic increasing (SI), right tail increasing (RTI) and positive quadrant dependence (PQD) properties. For more details, see  \cite{nelsen2007introduction}, \cite{karlin1968total}, and \cite{joe1997multivariate}. 

\par \cite{holland1987dependence} proved that the function $g(x,y)$ has $TP_2$ property if 
$$\zeta_g(x,y)=\frac{\partial^2 \ln g(x,y)}{\partial x\partial y} \geq 0, \quad \forall (x,y) \in \mathbb{R}^2.$$
Now, using the result of \cite{holland1987dependence}, we will prove the following result. 
\begin{proposition}
	The copula density (\ref{copulad}) has $TP_2$ property if $\delta\geq0$ and $\alpha\leq2$.
\end{proposition}
\begin{proof}
We have,
\begin{align*}
	\zeta_c(u,v)=&\dfrac{\partial^2 \ln c(u,v)}{\partial x\partial y}\\
	=&\dfrac{\delta\alpha^2\left[2-\alpha(1-2u)\right]\left[2-\alpha(1-2v)\right]\exp\left\{\alpha\left((u-u^2)+(v-v^2)\right)\right\}}{\left(1+\delta \alpha^2\left(1-2u\right)\left(1-2v\right)\exp\left\{\alpha\left((u-u^2)+(v-v^2)\right)\right\}\right)^2}\geq0,
\end{align*}
for every $\delta\geq0$ and $\alpha\leq2$.
\end{proof}
\subsection{Tail Dependence Coefficients}
The tail dependence coefficients measure the level of dependency among the random variables in the upper-right quadrant and in the lower-left quadrant of $\mathcal{I}^2$. In terms of copula, the upper and lower tail dependence coefficients, denoted by $\lambda_U$ and $\lambda_L$ respectively, are given by
\begin{equation}\label{tail}
\lambda_L=\lim_{u\rightarrow0^+}\frac{C(u,u)}{u}, \ \
\lambda_U=\lim_{u\rightarrow1^-}\frac{1-2u+C(u,u)}{1-u}.
\end{equation}
It is known that  $0\leq\lambda_L\leq1$ and $0\leq\lambda_U\leq1$. If $\lambda_L \in (0, 1]$, we say the copula $C(u,v)$ has lower tail dependence, and if $\lambda_L=0$, we say $C(u,v)$ has no lower tail dependence. Similar interpretation can be made for $\lambda_U$ (see \cite{nelsen2007introduction}, p. 214). Now, we will prove the following result. 
\begin{proposition}
The copula (\ref{Copula_prop}) has no tail dependence.
\end{proposition}
\begin{proof}
Using Eq. (\ref{tail}), we have
\begin{align*}
	\displaystyle\lambda_L&= \displaystyle\lim_{u\rightarrow0^+}\frac{u^2+\delta\left(1-e^{\alpha(u-u^2)}\right)^2}{u}\\
	&=\displaystyle\delta\lim_{u\rightarrow0^+}\frac{1}{u}\left(1-\sum_{n=0}^{\infty}\frac{\alpha^n u^n (1-u)^n}{n!}\right)^2\\
	&=\displaystyle\delta\lim_{u\rightarrow0^+}\alpha^2 u (1-u)^2\left(\sum_{n=2}^{\infty}\frac{\alpha^{n-1} u^{n-1} (1-u)^{n-1}}{n!}\right)^2 \\
&=0.
\end{align*}
\begin{align*}
	\displaystyle\lambda_U
	&=\displaystyle\lim_{u\rightarrow1^-}\frac{1-2u+u^2+\delta\left(1-e^{\alpha(u-u^2)}\right)^2}{1-u}\\
	&=\displaystyle\lim_{u\rightarrow1^-}\frac{(1-u)^2+\delta\left(1-\sum_{n=0}^{\infty}\frac{\alpha^n u^n (1-u)^n}{n!}\right)^2}{(1-u)}\\
	&=\displaystyle\delta\lim_{u\rightarrow1^-}\alpha^2 u^2 (1-u)\left(\sum_{n=2}^{\infty}\frac{\alpha^{n-1} u^{n-1} (1-u)^{n-1}}{n!}\right)^2 \\
	&=0.
\end{align*}
\end{proof} In the next section, we will develop a bivariate Rayleigh distribution as an application of the proposed copula (\ref{Copula_prop}). We will study some statistical properties of new bivariate Rayleigh distribution, and a real data analysis involving the new distribution is also presented.

\section{A New Bivariate Rayleigh distribution}\label{4}
\begin{figure}[ht]
	\centering
	\subfloat[]{\includegraphics[height=5cm, width=8cm]{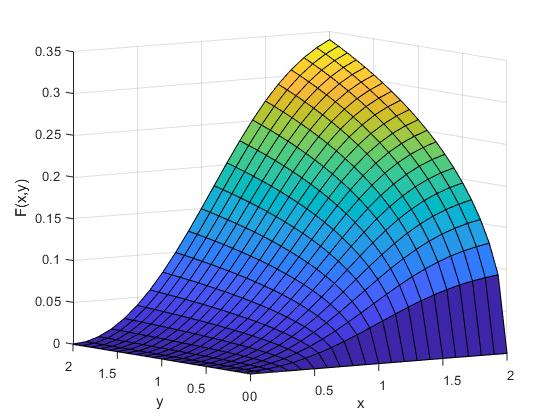}}            
	\subfloat[]{\includegraphics[height=5cm, width=10cm]{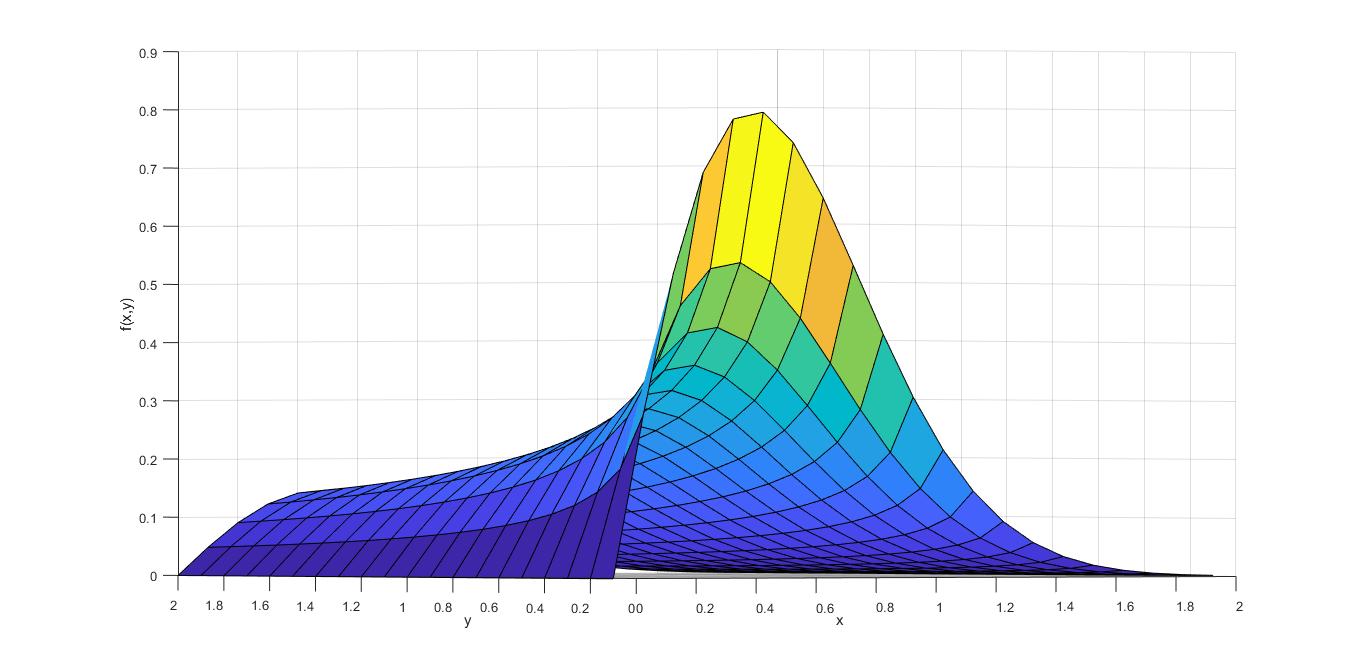}}  
	\caption{Surface plots of $F(x,y)$ and $f(x,y)$ of the BRD distribution for $\lambda_{1}=3$, $\lambda_{2}=2$, $\delta=0.5$, $\alpha=3.8$. 
	 }\label{cdfplot}
\end{figure}

Rayleigh distribution is one of the most popular models in medical sciences, engineering, particle physics and economics. A random variable $X$ follows Rayleigh distribution with parameter $\lambda$, denoted by $Rayleigh(\lambda)$, if its cumulative distribution function (CDF) is given by
$F(x;\lambda)=1-e^{-x^2/2\lambda^2},\ x>0,\lambda>0$, and corresponding probability density function (PDF) is given by $f(x;\lambda)=\frac{x}{\lambda^2} e^{-x^2/2\lambda^2}, \ x>0,\lambda>0$. Let $X$ and $Y$ be two random variables having $Rayleigh(\lambda_1)$ and $Rayleigh(\lambda_2)$ distributions respectively, and the dependence between $X$ and $Y$ is modelled by the copula (\ref{Copula_prop}). Then, the join distribution function of $X$ and $Y$ is given by
\begin{align}\label{BRD}
 F(x,y;\Theta)=&\left(1-e^{-x^2/2\lambda_1^2}-e^{-y^2/2\lambda_2^2}+e^{-\left(x^2/2\lambda_1^2+y^2/2\lambda_2^2\right)}\right) +\delta\left(1-e^{\alpha\left(e^{-x^2/2\lambda_1^2}-e^{-x^2/\lambda_1^2}\right)}\right)\nonumber\\&\hspace{1cm}\left(1-e^{\alpha\left(e^{-y^2/2\lambda_2^2}-e^{-y^2/\lambda_2^2}\right)}\right),
\end{align}
where $x>0,y>0,\lambda_1>0,\lambda_2>0,\alpha\in \mathbb{R},\lvert\delta\rvert\leq \delta ^{\star}(\alpha)$ and $\Theta=\left(\lambda_1,\lambda_2,\alpha,\delta\right)$. A non-negative random vector $(X,Y)$ is said to follow bivariate Rayleigh distribution with parameters $\lambda_1,\lambda_2,\alpha$ and $\delta$, if its joint CDF is given by Eq. (\ref{BRD}) and is denoted by BRD$(\lambda_1,\lambda_2,\alpha,\delta)$. The corresponding joint density function is given by
\begin{align}\label{BRD1}
f(x,y;\Theta)=&\left( \frac{xy}{\lambda_1^2\lambda_2^2} e^{-\left(x^2/2\lambda_1^2+y^2/2\lambda_2^2\right)}\right)\displaystyle \left[1+\delta \alpha^2\left(2e^{-x^2/2\lambda_1^2}-1\right)\left(2e^{-y^2/2\lambda_2^2}-1\right)\right.\nonumber\\
	&\quad\quad\left.
	\left(\exp\left\{\alpha\left(e^{-x^2/2\lambda_1^2}-e^{-x^2/\lambda_1^2}+e^{-y^2/2\lambda_2^2}-e^{-y^2/\lambda_2^2}\right)\right\}\right)\right],
\end{align}
where $x>0,y>0,\lambda_1>0,\lambda_2>0,\alpha\in \mathbb{R},\lvert\delta\rvert\leq \delta ^{\star}(\alpha)$ and $\Theta=\left(\lambda_1,\lambda_2,\alpha,\delta\right)$. Surface plots of joint CDF (\ref{BRD}) and joint density function (\ref{BRD1}) of the BRD family are shown in Figure \ref{cdfplot}. These figures are constructed using Matlab R2021b. 
 Now, we will provide expressions for conditional distribution and product moments.
\begin{proposition}
	Let $(X,Y)\sim BRD(\lambda_1,\lambda_2,\alpha,\delta)$. Then
	\begin{enumerate}[(i)]
		\item $X\sim Rayleigh(\lambda_1)$ and $Y\sim Rayleigh(\lambda_2),$
		\item the conditional density function of $X$ given $Y=y$ is
		\begin{align*}
			f(x|y)=&\left( \frac{x}{\lambda_1^2} e^{-x^2/2\lambda_1^2}\right)
			\displaystyle \left[1+\delta \alpha^2\left(2e^{-x^2/2\lambda_1^2}-1\right)\left(2e^{-y^2/2\lambda_2^2}-1\right)\right.\nonumber\\
			&\left.\hspace{0.5cm}\left(e^{\alpha\left(e^{-x^2/2\lambda_1^2}-e^{-x^2/\lambda_1^2}+e^{-y^2/2\lambda_2^2}-e^{-y^2/\lambda_2^2}\right)}\right)\right],
		\end{align*}
		\item the conditional distribution function of $X$ given $Y=y$ is
		\begin{align*}
			F(x|y)=&\left(1-e^{-x^2/2\lambda_1^2}\right)+\delta \alpha\left(e^{\alpha\left(e^{-x^2/2\lambda_1^2}-e^{-x^2/\lambda_1^2}\right)}-1\right)\left(2e^{-y^2/2\lambda_2^2}-1\right)\left(e^{\alpha\left(e^{-y^2/2\lambda_2^2}-e^{-y^2/\lambda_2^2}\right)}\right),
		\end{align*}
		where $x>0,y>0,\lambda_1>0,\lambda_2>0,\alpha\neq0,\lvert\delta\rvert\leq \delta ^{\star}(\alpha)$ and $\Theta=\left(\lambda_1,\lambda_2,\alpha,\delta\right)$.
	\end{enumerate}
\end{proposition}
\begin{proposition}\label{prdmnt}
	Let $(X,Y)\sim BRD(\lambda_1,\lambda_2,\alpha,\delta)$. Then $(r,s)$-th order product moments can be expressed as
	\begin{align*}
		E(X^rY^s)=&\lambda_1^r\lambda_2^s2^{(r+s)/2}\Gamma\left(1+r/2\right)\Gamma\left(1+s/2\right)\\
		&\hspace{.5cm}\left[1+\delta\alpha^2\left(\sum_{k=0}^{\infty}\frac{\alpha^k}{k!}\sum_{t=0}^{k}\frac{(-1)^t}{\lambda_1^2}\binom{k}{t}\left(2\left(k+t+2\right)^{-\frac{r+2}{2}}-\left(k+t+1\right)^{-\frac{r+2}{2}}\right)\right)\right.\\
		&\hspace{.5cm}\left.\left(\sum_{k=0}^{\infty}\frac{\alpha^k}{k!}\sum_{t=0}^{k}\frac{(-1)^t}{\lambda_2^2}\binom{k}{t}\left(2\left(k+t+2\right)^{-\frac{s+2}{2}}-\left(k+t+1\right)^{-\frac{s+2}{2}}\right)\right)\right],
	\end{align*}
	where $\Gamma(t)$ denotes the usual gamma function. 
\end{proposition}
\subsection{Real Data Application}
We consider the UEFA Champions League data set from 2004 to 2006, reported in  \cite{meintanis2007test}. In this data set, $X$ and $Y$ represent the time (in minutes) of the first goal scored by Team-A and Team-B, respectively. To check whether the marginal distributions of $X$ and $Y$ support the Rayleigh distribution, we perform Kolmogorov-Smirnov (KS) one-sample test. The KS-test's results suggest that $X$ supports Rayleigh distribution with parameter $\hat{\lambda_1}=32.14599$ (p-value=0.934 and KS statistic value=0.088515). Similarly, $Y$ also supports Rayleigh distribution with parameter $\hat{\lambda_2}=28.172255$ (p-value=0.07727 and KS statistic value=0.20968). We fit the proposed bivariate Rayleigh distribution, and the results are shown in Table \ref{bivariateMLE}. We compare the new BRD model with Marshall Olkin's bivariate exponential distribution (BMOED) by \cite{meintanis2007test}, bivariate generalized exponential distribution (BGED) by \cite{mirhosseini2015new}, and bivariate generalized Rayleigh distribution (BGRD) proposed by \cite{pathak2020bivariate}. We use the log-likelihood (LL) function, Akaike Information Criteria (AIC) and Bayesian Information Criteria (BIC) as the comparison criteria. The formulas for AIC and BIC are given by 
\[\text{AIC}=2k-2\ln L, \quad \text{and} \quad \text{BIC}=k\ln n-2\ln L,\]
where $k$ is the number parameters in the model, $n$ is the sample size and $L$ is the maximum value of the likelihood function. From Table \ref{bivariateMLE}, it is clear that bivariate Rayleigh distribution provides a better fit over BGED, BMOED and BGRD for the UEFA champions league data set. 
\begin{table}[h]
	
	\caption{ML estimates, LL, AIC, and BIC values for the bivarite distributions using  UEFA Champion's League data set.} 
	\centering 
	\scalebox{.85}{\begin{tabular}{c c c c c} 
			\hline 
			Bivariate Distribution & ML Estimates & LL  & AIC & BIC \\ [0.5ex] 
			\hline 
			BGED & $\hat{\alpha}_1=0.0244, \hat{\alpha}_2=0.0304, \hat{\theta}=0.999$  & -340.5234 & 687.0468 & 691.8795 \\[0.5ex]
			BMOED & $\hat{\lambda}_1=0.012, \hat{\lambda}_2=0.014, \hat{\lambda}_3=0.022$ & -339.006 & 684.012 & 688.8448 \\[0.5ex]
			BGRD & $\hat{b}_1=0.000530, \hat{b}_2=0.000836, \hat{\theta}=0.40331$ & -331.879 & 664.589 & 672.6436 \\[0.5ex]
			BRD & $\hat{\lambda}_1=33.39429, \hat{\lambda}_2=28.08949, \hat{\delta}=10.39829, \hat{\alpha}=0.2871858 $ & -327.256 & 664.512 & 668.9557 \\[0.5ex]
			\hline 
	\end{tabular}}
	\label{bivariateMLE} 
\end{table}
\section{Conclusion and Future Direction}
This paper proposes a new bivariate symmetric copula exhibiting positive and negative dependence. The main features of the copula are: (i) it has a simple mathematical structure, (ii) it has a wider dependence range when compared to FGM copula and its generalizations, and (iii) there is no lower and upper tail dependence. Using the proposed copula, we developed a new bivariate Rayleigh distribution (BRD) and discussed some statistical properties. The proposed bivariate model provides a better fit for a  real data set. Since we considered only the symmetric version of the bivariate copula, the asymmetric version is still an open problem for new researchers. 

\section*{Acknowledgments}
The first author is thankful to the Indian Institute of Technology Indore for providing financial support for pursuing research.
\subsection*{Conflict of interest}
The authors declare no potential conflict of interests.
\bibliographystyle{chicago}
\bibliography{references}



\end{document}